\documentclass[12pt]{article}
\usepackage{amssymb, amsmath}
\begin{document}
\newtheorem{theorem}{Theorem}
\newtheorem{lemma}{Lemma}[section]
\newtheorem{definition}{Definition}
\newtheorem{pro}{Proposition}
\newtheorem{cor}{Corollary}
\newcommand{\n}{\nonumber}
\newcommand{\tv}{\tilde{v}}
\newcommand{\tw}{\tilde{\omega}}
\renewcommand{\t}{\theta}
\newcommand{\w}{\omega}
\newcommand{\e}{\varepsilon}
\renewcommand{\a}{\alpha}
\renewcommand{\l}{\lambda}
\newcommand{\vare}{\varepsilon}
\newcommand{\s}{\sigma}
\renewcommand{\o}{\omega}
\newcommand{\bb}{\begin{equation}}
\newcommand{\ee}{\end{equation}}
\newcommand{\bq}{\begin{eqnarray}}
\newcommand{\eq}{\end{eqnarray}}
\newcommand{\bqn}{\begin{eqnarray*}}
\newcommand{\eqn}{\end{eqnarray*}}
\title{A note on  `Nonexistence of  self-similar singularities for
the 3D incompressible Euler equations'}
\author{Dongho Chae\thanks{The work was supported
partially by the KOSEF Grant no. R01-2005-000-10077-0.} \\
Department of Mathematics\\
              Sungkyunkwan University\\
              Suwon 440-746, Korea\\
  e-mail: {\it chae@skku.edu}}
 \date{}
\maketitle
\begin{abstract}
In this brief note we show that the author's previous result in
\cite{cha} on the nonexistence of  self-similar singularities for
the 3D incompressible Euler equations implies actually the
nonexistence of `locally self-similar' singular solution, which has
been sought by many physicists.
\end{abstract}

\section*{Nonexistence of locally self-similar solution}
We are  concerned here on the following Euler equations for the
homogeneous incompressible fluid flows on  $\Bbb R^3$.
 \[
\mathrm{ (E)}
 \left\{ \aligned
 &\frac{\partial v}{\partial t} +(v\cdot \nabla )v =-\nabla p ,
 \quad (x,t)\in {\Bbb R^3}\times (0, \infty) \\
 &\quad \textrm{div }\, v =0 , \quad (x,t)\in {\Bbb R^3}\times (0,
 \infty)\\
  &v(x,0)=v_0 (x), \quad x\in \Bbb R^3
  \endaligned
  \right.
  \]
where $v=(v_1, v_2, v_3 )$, $v_j =v_j (x, t)$, $j=1,2,3$, is the
velocity of the flow, $p=p(x,t)$ is the scalar pressure, and $v_0
$ is the given initial velocity, satisfying div $v_0 =0$.

The  system (E) has the scaling property that
  if $(v, p)$ is a
solution of the system (E), then for any $\lambda >0$ and $\alpha
\in \Bbb R $ the functions
 \bb
 \label{self}
  v^{\lambda, \alpha}(x,t)=\lambda ^\alpha v (\lambda x, \l^{\a +1}
  t),\quad p^{\l, \a}(x,t)=\l^{2\a}p(\l x, \l^{\a+1} t )
  \ee
  are also solutions of (E) with the initial data
  $ v^{\lambda, \alpha}_0(x)=\lambda ^\alpha v_0
   (\lambda x)$.

This  scaling
  property leads to the  following definition of self-similar blowing up
  solution  of (E):
  \begin{definition}
 A solution $v(x,t)$ of the solution to (E) is called a self-similar
 blowing up solution if there exist $\a> -1$, $T_* >0$ and a solenoidal vector field
 $V$ defined on $\Bbb R^3$ such that
  \bq
  \label{vel}
 v(x, t)=\frac{1}{(T_*-t)^{\frac{\a}{\a+1}}}
V\left(\frac{x}{(T_*-t)^{\frac{1}{\a+1}}}\right)\quad \forall
 (x,t)\in \Bbb R^3 \times (-\infty, T_* ).
 \eq
 \end{definition}
The above definition is apparently `global' in the sense that the
self-similar representation of solution in (\ref{vel}) should hold
for all space-time points in $\Bbb R^3 \times (-\infty, T_* )$. On
the other hand, many physicists have been trying to seek a
`locally self-similar' solution of the 3D Euler equations(see e.g.
\cite{pel} and the references therein). We formulate the precise
definition of this in the following.
\begin{definition}
 A solution $v(x,t)$ of the solution to (E) is called a locally self-similar
 blowing up solution near a space-time point
 $(x_*, T_*)\in \Bbb R^3 \times (-\infty, +\infty)$
 if there exist $r>0$, $\a> -1$ and a solenoidal vector field
 $V$ defined on $\Bbb R^3$ such that the
 representation
\bq
  \label{vela}
 v(x, t)=\frac{1}{(T_*-t)^{\frac{\a}{\a+1}}}
V\left(\frac{x-x_*}{(T_*-t)^{\frac{1}{\a+1}}}\right) \,\,
 \forall (x,t)\in  B(x_*, r)\times (T_*-r^{\a+1}, T_* )
 \eq
  holds true, where
 $ B(x_*, r)=\{ x\in \Bbb R^3 \, |\, |x-x_* |<r \}$.
 \end{definition}

We have the following relation between the two notions of
self-similar blowing  up solutions.
\begin{theorem}
The nonexistence of  self-similar solution of the 3D Euler equations
in the sense of Definition 1 implies the nonexistence of  locally
self-similar solution in the sense of Definition 2.
\end{theorem}

Combining Theorem 1 with the main theorem in \cite{cha}(Theorem
1.1), we have the following corollary.

\begin{cor} Suppose there exists a locally self-similar blowing up
solution of the 3D Euler equations in  the form (\ref{vela}). If
there exists $p_1>0$ such that $\Omega =\mathrm{curl}\, V \in
L^p(\Bbb R^3 )$ for all $p\in (0, p_1)$, then necessarily $\Omega
=0$. In other words, there exists no nontrivial locally self-similar
solution to the 3D Euler equation if the vorticity $\Omega$
satisfies the integrability condition above.
\end{cor}

\noindent{\bf Proof of Theorem 1.} We assume there exists a locally
self-similar blowing up solution $v(x,t)$ in the sense of Definition
2. The proof of the theorem follows if we prove the existence of
self-similar blowing up solution in the sense of Definition 1 based
on that assumption.
 By translation in
space-time variables, we can rewrite the velocity
 in (\ref{vela}) as
 \bb
  \label{velb}
 v(x, t)=\frac{1}{t^{\frac{\a}{\a+1}}}
V\left(\frac{x}{t^{\frac{1}{\a+1}}}\right)\quad \mbox{for}\quad
(x,-t)\in B(0,r)\times (-r^{\a+1}, 0).
 \ee
  We observe that, under the scaling
transform (\ref{self}), we have the invariance of the
representation,
$$
v(x, t) \mapsto v^{\lambda, \alpha} (x,t)=\lambda ^{\alpha}
v(\lambda x, \l^{\a +1}
  t)=\frac{1}{t^{\frac{\a}{\a+1}}} V
\left(\frac{x}{t^{\frac{1}{\a+1}}}\right) (=v (x,t)),
$$
while the region of space-time, where the self-similar form of
solution is valid, transforms according to
$$
B(0,r)\times (-r^{\a+1}, 0) \mapsto B(0,r/\lambda) \times
\left(-(r/\lambda)^{\a+1}, 0\right) .$$

We set $\l =1/n$, and define the sequence of locally self-similar
solutions $\{ v^n (x,t)\}$  by $v^n (x,t):=v^{\frac{1}{n}, \a}(x,t)$
with $v^1 (x,t)=v(x,t)$. In the above we find that
$$
v^n (x,t)=\frac{1}{t^{\frac{\a}{\a+1}}} V
\left(\frac{x}{t^{\frac{1}{\a+1}}}\right) \quad \mbox{for} \quad
 (x,-t)\in B(0, nr)\times ( -(nr)^{\a+1},0),
$$
and each $v^n (x,t)$ is a solution of the Euler equations for all
$(x,t)\in \Bbb R^3 \times (-\infty, 0)$. Let us define $v^\infty
(x,t)$ by
$$v^\infty (x,t)=\frac{1}{t^{\frac{\a}{\a+1}}} V
\left(\frac{x}{t^{\frac{1}{\a+1}}}\right)\quad \mbox{for}\quad
 (x,-t)\in \Bbb R^3 \times ( -\infty,0).
 $$
 Given a compact set $K \subset \Bbb R^3\times
(-\infty ,0) $, we observe that $v^n \to v^\infty$ as $n\to \infty$
on $ K$ in any strong sense of convergence. Indeed, for sufficiently
large $N=N (K)$,  $v^n (x,t)\equiv v^\infty (x,t)$ for $(x,t)\in K$,
if $n>N$. Hence, we find that $v^\infty(x,t)$ is a solution of the
Euler equations for all $(x,t)\in \Bbb R^3 \times (-\infty , 0)$,
which is a self-similar blowing up solution, after translation in
time, in the sense of Definition 1.
 $\square$\\

 We note that the above proof obviously works also for the
 self-similar solutions of the other equations considered in \cite{cha} and
  Leray's self-similar solutions of the Navier-Stokes equations(see \cite{ler} for the problem,
   and \cite{nec} for the rule-out of the `global' self-similar
 solution).
 \[ \mbox{\bf Acknowledgements}\]
 The author would like to thank to Prof. K. Ohkitani and
Dr. T. Matsumoto  at Kyoto University for useful discussions, and
informing him of the reference\cite{pel}.

\end{document}